\documentclass[12pt]{amsart}
\usepackage{amsthm,amsmath,amssymb}
\usepackage{addfont}
\usepackage{mathrsfs}
\usepackage[T1]{fontenc}
\numberwithin{equation}{section}

\def\ca{{\mathcal A}}

\def\cd{{\mathcal D}}

\def\cq{{\mathcal Q}}
\def\car{{\mathcal R}}
\def\cas{{\mathcal S}}

\def\bc{{\mathbb C}}

\def\bn{{\mathbb N}}

\def\br{{\mathbb R}}

\def\bt{{\mathbb T}}

\def\eps{\varepsilon}

\def\d{\delta}

\def\g{\gamma}

\def\l{\lambda}

\def\s{\sigma}

\def\f{\varphi}

\theoremstyle{plain}
\newtheorem{lemma}{Lemma}[section]
\newtheorem{proposition}[lemma]{Proposition}
\newtheorem{theorem}[lemma]{Theorem}
\newtheorem{corollary}[lemma]{Corollary}
\theoremstyle{definition}

\newtheorem{definition}[lemma]{Definition}

\begin{document}

\title[Column bounded matrices and Grothendieck's inequalities]{\textsc{ Column bounded matrices and Grothendieck's inequalities  }}

\author[E.~Christensen]{Erik Christensen}
\address{\hskip-\parindent
Erik Christensen, Mathematics Institute, University of Copenhagen, Copenhagen, Denmark.}
\email{echris@math.ku.dk}
\date{\today}

\begin{abstract} It follows from Grothendieck's {\em little inequality} that to any complex $m \times n $ matrix $X$ of column norm at most 1, and an $0 <\eps <1,$ there exist a natural number $l, $ a matrix $C$ in $M_{(m,l)}(\bc)$  with $(1-\eps)^2I_m \leq CC^* \leq (4/\pi) (1 + \eps)^2 I_m $  and an $l \times n$ matrix $Z$ with entries in the torus $\bt,$ such that $X=  l^{-(1/2)}CZ.$    
Grothendieck's {\em little   inequality } with the constant $k_G^\bc = 4/\pi $ follows  from this factorization. Grothendieck's {\em big inequality} may be reproduced with the estimate $K_G^\bc \leq k_G^\bc/(2- k_G^\bc).$  
 \end{abstract}

\maketitle

\section{Introduction}

Grothendieck's work on tensor products of Banach spaces \cite{Gr} has influenced mathematics in several ways, some of which are very surprising. This is described in Pisier's survey article \cite{Pi3}.  Here we will focus on the inequality, which is named   Grothendieck's {\em little inequality} in the setting of complex $m \times n$ matrices. That inequality is established in Section 3 of \cite{Gr} and is presented as a factorization result of a bounded operator from a Hilbert space to an $L^1(\mu)$ space. In this article all the vector spaces we deal with are complex and  finite dimensional, but nevertheless we will equip the spaces with norms, such that we are still working with the same concepts as Grothendieck studied. 

The results we present here are based on Grothendieck's insights and  on  our recent articles \cite{C1}, \cite{C3}, where we showed that the theory of operator spaces and completely bounded maps provides a set up, which fits very well - in our opinion - to the existing results related to Grothendieck's inequalities. When we worked on these articles we were unaware of the article \cite{PV} from 2016,  where C. Palazuelos and T.  Vidick in a survey article demonstrates that the language of operator spaces and completely bounded norms fits nicely in to the study of questions related to {\em Bell's inequality}  in quantum mechanics. This inequality has its offspring in the {\em Einstein, Podolsky, Rosen paradox, } which mathematically is based on the use of {\em tensor products } in the  quantum mechanical description of the  coupling of two  systems.
The Theorem 11.12 in  \cite{AS}, gives a concrete translation of Grothendieck's {\em big inequality } over the real numbers  into a statement on the following  relation between two concepts from quantum mechanics. The  set of local correlation matrices,  LC$_{m,n},$ and the set of  quantum correlation matrices QC$_{m,n} $ satisfy the relation QC$_{m,n} \subseteq K_G^\br ($LC$_{m,n}).$  This connection between Bell inequalities and Grothendieck's big inequality was first described by Tsirelson \cite{Ts}.     In this article we look at complex vector spaces only, and we consider Minkowski norms associated  to convex sets rather than convex subsets of a vector space. In our set up  the analogy to the  statement on correlation matrices may be found in equation (\ref{ScGr}).

  The set up in this article is basically a study of the tensor product $\bc^m \otimes \bc^n,$  with the extra condition that $\bc^m$ will be equipped with the norm  $\|.\|_p $ and $\bc^n$ will be equipped with the norm $\|.\|_r$ for  $p, r $ in $\{1, 2, \infty\}.$   We will express this study as a study of the complex $m \times n $ matrices, which we denote $M_{(m,n)}(\bc).$  
This linear space  of matrices is in a canonical way isomorphic to the algebraic tensor product $\bc^m \otimes \bc^n $ where the isomorphism  is described via the canonical basis $(\d_i)_{\{1 \leq i \leq m\}}$ for  $\bc^m$ and $(\g_j)_{\{1 \leq j\leq n\}},$ for $\bc^n$ and the matrix units $\{e_{(i,j)} : 1 \leq i \leq m,\, \, 1 \leq j \leq n\}$ for $M_{(m,n)}(\bc)$ by  the linear  map $\f : \bc^m \otimes \bc^n \to M_{(m,n)}(\bc ) $ which satisfies $\f(\d_i \otimes \g_j ):= e_{(i,j)}.$
In several spots we will use that the image $\f(a \otimes b) $ is the rank one matrix with entries $\f(a \otimes b)_{(i,j)}= a_ib_j$ and also use, that this matrix is a product of a one column matrix $a_| := (a_1, \dots, a_m)$ and a one row matrix $b_{-} := (b_1, \dots, b_n),$ so $\f(a \otimes b ) = a_|b_{-}.$

 Given a couple of normed spaces   such as $(\bc^m, \|.\|_p)$ and $(\bc^n, \|.\|_r), $ we recall that Schatten \cite{Scha} has introduced the concept named a cross norm on the tensor product $(\bc^m, \|.\|_p) \otimes (\bc^n, \|.\|_r),$ and we recall that a norm say $|||.|||$  on the tensor product of the normed spaces is called a cross norm on this tensor product of normed spaces if it satisfies $$ \forall \eta \in \bc^m \, \forall\xi \in \bc^n : \quad |||(\eta \otimes \xi)|||  = \|\eta\|_p \|\xi\|_r.$$ 
Schatten proved that there is a minimal and a maximal cross norm. Today the minimal cross norm  is called the {\em injective cross norm } and - in the case above - it is denoted  $\|.\|_{\vee(p,r)}. $ The maximal cross norm is called the {\em projective cross norm } and - in the case above - it is denoted $\|.\|_{\wedge(p,r)}.$  Via the isomorphism $\f$ we will think of theses norms as norms on $M_{(m,n)}(\bc).$ 

 There are other well known norms on $M_{(m,n)}(\bc),$ and amongst them we will right now  mention the {\em operator norm, } which we denote $\|X\|_\infty = \|X\|_{\vee(2,2)},$ and the {\em Hilbert Schmidt norm,}  which we denote $\|X\|_2,$ is also a $(2,2)$ cross norm.   
 
 We will remind you on the Schur product on the complex  matrices. Given two complex matrices $X = (X_{(i,j)}) $ and $A= (A_{(i,j)} )$ in $M_{(m,n)} (\bc),$ then we define their Schur product $X \circ A$ to be the matrix in $M_{(m,n)}(\bc), $ which is the entry-wise product of the matrices,  $(X \circ A)_{(i,j)} := X_{(i,j)}A_{(i,j)}.$ 
 The norm of $X$ when it acts as a Schur multiplier on $(M_{(m,n)}(\bc), \|.\|_\infty)$ is denoted $\|X\|_S, $ and the norm $\|.\|_S$  is an $(\infty, \infty)$ cross norm, which has the property that its completely bounded version $\|.\|_{cbS}$ equals  $\|.\|_S,$ \cite{Sm}, .
 
  It is possible, \cite{Pi3} Theorem 3.2, to formulate Grothendieck's {\em big inequality} in the language of cross norms on tensor products as  \begin{align} \label{ScGr} &\text{ There exists a minimal positive real } K_G^\bc \text{  such that } \\ \notag &\forall m,n  \in \bn\, \forall X \in M_{(m,n)}(\bc): \quad \|X\|_{\wedge(\infty, \infty)} \leq K_G^\bc \|X\|_S. \end{align}  
 Grothendieck's {\em big inequality}  is most often described as a property for bilinear forms  on the product of two abelian C*-algebras. Here we  study finite dimensional C*-algebras,  and then such a bilinear form may be viewed as a linear form on $(\bc^m, \|.\|_\infty )\otimes (\bc^n, \|.\|_\infty),$ When a complex $m \times n$ matrix $X$ acts as a kernel for such a bilinear form we denote the norm of the bilinear  form by $\|X\|_B,$ and we note that $\|X\|_B = \|X\|_{\vee(1,1)} .$ This bilinear form  also has  a completely bounded norm, as described in \cite{CS}, \cite{Pa}, \cite{C3}, which  we denote by  $\|X\|_{cbB},$  and we recall that $\|.\|_{cbB}$ is a $(1,1)$ cross norm, \cite{C1} Theorem 2.1.   We will return to this after a short reminder on the inner product $\langle X, Y\rangle $ on $M_{(m,n)}(\bc).$ This inner product is defined, via that trace Tr$_n$ on $M_n(\bc)$ which has the property  Tr$_n(I_n) = n,$ by 
$$\langle X, Y\rangle := \mathrm{Tr}_n (Y^*X).$$
We showed in \cite{C1},  that via this sesqui-linear form, the norms $\|.\|_S $ and $\|.\|_{cbB}$ are conjugate dual to each other. Schatten shows that the minimal and the maximal cross norms are dual to each other, so the equation (\ref{ScGr}) may be dualized into the equivalent form 
 \begin{align} \label{gF} &\text{ There exists a minimal positive real } K_G^\bc \text{ such that }\\ \notag 
\forall m,n & \in \bn\,  \forall X \in M_{(m,n)}(\bc): \quad \|X\|_{cbB} \leq K_G^\bc \|X\|_{\vee(1,1)} = K_G^\bc \|X\|_B. 
\end{align}   
 From this inequality and the theorems on the structure of completely bounded linear or multilinear maps \cite{Pa}, we can obtain   some new and some  well known factoirzation results as described in \cite{C3} and \cite{C4}.

Grothendieck's {\em little inequality} may also be expressed as a relation between a norm and its completely bounded version. 
To  formulate that we recall the norms $\|X\|_F$ and $\|X\|_{cbF} $  from   \cite{C3},  which we defined on $M_{(m,n)}(\bc).$ The norm  $\|X\|_F$  is defined as the norm of the linear operator $F_X$ with the matrix $X$ acting as an operator from $(\bc^n, \|.\|_\infty)$ to $(\bc^m, \|.\|_2),$ so $\|X\|_F = \|X\|_{\vee(2,1)}$  and then $\|X\|_{cbF}$ is defined as the completely bounded norm of $F_X,$ and it follows that $\|.\|_{cbF} $ is a $(2, 1)$ cross norm, \cite{C1} Theorem 2.1.   It follows from \cite{Pi3} Section 5 and  Theorem 3.11 of  \cite{C3} that Grothendiek's {\em little  inequality} may be formulated as: 
\begin{align} \label{lG} & \text{{\em The constant $k_G^\bc = \frac{4}{\pi} $  is the minimal positive real such that}}\\
\notag
\forall m,n \in &\bn,\, \forall X \in M_{(m,n)}(\bc): \,\,\|X\|_{cbF} \leq \sqrt{k_G^\bc}\|X\|_F = \sqrt{k_G^\bc}\|X\|_{\vee(2,1)}. 
\end{align} 
For a complex $m \times n$ matrix $X$ we defined in  \cite{C3}, item (iv) on the first page of \cite{C3},  a bilinear operator $T_X: (\bc^m, \|.\|_\infty) \times M_{(m,n)}(\bc) \to \bc^n.$  It turned out, \cite{C3} equation (3.6), that the  norm $\|X\|_T:= \|X_T\|$  is  conjugate dual to $\|.\|_{cbF}$ under the duality implied by the inner product $\langle X, Y\rangle,$ and it follows that $\|.\|_T$ is a $(2, \infty)$ cross norm, which satisfies $\|.\|_T = \|.\|_{cbT}, $ \cite{C1} Theorem 2.1.   By duality, we can then obtain an equivalent formulation of the {\em little inequality} as 
\begin{align}  \label{gT} & \text{ The  constant }k_G^\bc = \frac{4}{\pi} \text{ is the minimal positive real such that } \\ \notag
\forall m,n & \in \bn \,\forall X \in M_{(m,n)}(\bc):\quad \|X\|_{\wedge(2, \infty)} \leq \sqrt{k_G^\bc}\|X\|_T.
\end{align} 

\medskip
The results on column bounded matrices, which we present here, are based on our studies in \cite{C1}, were  we studied some  relations between a couple of compact convex sets of positive complex $n \times n$ matrices $\cq_n$ and $\car_n. $ These sets  are defined below, and we have included yet another compact convex set $\cd_n $ of positive matrices.  

\begin{definition}
Let $n $ be a natural  number, the compact convex subsets $\cd_n, \cq_n, $ and $\car_n$   of the positive matrices in $M_n(\bc)$ are defined by 
\begin{itemize}
\item[(i)] $\cd_n := \{D \in M_n(\bc): D \geq 0 \text{ and diag}(D) \leq I_n\}.$ 
\item[(ii)] $\cq_n := \{Q \in M_n(\bc)\, : \, Q \geq 0 \text{ and } \mathrm{diag}(Q) = I_n\}.$ 
\item[(iii)] The set $\car_n$ is the closed convex hull of  the positive rank one matrices which have the unit $I_n$ as their main diagonal. 
\end{itemize}
\end{definition} 

Based on Grothendieck's little inequality - in the complex case - we showed in     
Proposition 5.3 of \cite{C1}  a slightly weaker result than the following proposition.

 \begin{proposition} \label{geo}  
The constant $k_G^\bc = 4/\pi$ is the smallest possible posi-tive real, such that for any natural number $n$ and any matrix $D$  in $\cd_n$ there exists a matrix   $R$ in $\car_n$ with $  D \leq   k_G^\bc R.$     
 \end{proposition}   
 
 \begin{proof}
 The proposition 5.3 of \cite{C1} is here  extended in  the way, that the original proposition is identical to the proposition above except that now $\cd_n$ has replaced $\cq_n$ in the original version.  To any matrix $D$ in $\cd_n$ you can find a positive diagonal $n \times n $ matrix $C$ such that $D+C $ is in $\cq_n$ and then the proposition follows. 
 \end{proof}
 
 This proposition is the basis for the results on column bounded matrices which we obtain, the reason being that for an $ r \times n$ matrix $X$ with column  norm at most 1, we clearly have that $X^*X $ is in $\cd_n$ and then Proposition \ref{geo} and  a well known elementary result in operator theory shows that $X$ may be written as $X=CR^{(1/2)}$ for an $R$ in $\car_n$ and a $C$ with $\|C\|_\infty \leq \sqrt{k_G^\bc} .$ This factorization may be written in a neater form, once the structure of $R$ is taken into account, and we present these results as Theorem \ref{ColBdFac} and Theorem \ref{ColBd2} below. 
These theorems are then consequences of Grothendieck's little inequality, and we will show that Grothendieck's little inequality may be obtained as a corollary to Theorem \ref{ColBdFac}. We present in Theorem \ref{ColBd2} a slightly weaker version of Theorem \ref{ColBdFac}, which has the advantage that the factorization we obtain in the latter case is simpler and easier to use in possible applications, we hope.  

The Theorem \ref{ColBd2} does not depend on the number of columns $n,$ so  if we expand the matrix $X$ to an $m \times (n+p)$ matrix with colun norm at most one and then apply Theorem \ref{ColBd2} to the expanded matrix,  we will still obtain a factorization of $X,$ but it will also give a factorization of the expanded matrix, both using the same matrix $C$ in the factorization.   
  We have not been able to make any real  use of this observation, but Theorem \ref{ManyCol} offers some consequences of this observation, and we will like to mention, that the latter theorem shows that it is always  possible to obtain a matrix $C$ in the factorization  which satisfies $$ (1- \eps)^2 I_m  \leq CC^* \leq (1+\eps )^2 k_G^\bc I_m,$$
for any $0 < \eps < 1.$ 

We end the article with a section that explicitly shows how both  of Grothendieck's inequalities  may be obtained as  consequences of the factorization result in  Theorem \ref{ColBdFac}. This application shows that $K_G^\bc \leq k_G^\bc/(2- k_G^\bc) < 1.752,$ an observation which also  was obtained in equation (4.6) of \cite{C1}. 
 \section{Factorizations of column bounded matrices}

We have used the word {\em column norm} without a formal definition because we think that the concept is well known, but it may not be so. Further we will like to use some notation of our own.

\begin{definition} Let $X$ be a complex $m \times n  $ matrix, and $i, j $ natural numbers in respectively 
$\{1, \dots, m\}$ and $\{1, \dots, n \}$
\begin{itemize}
\item[(i)] $\forall i: \, _iX$ denotes the $i$´th row of $X.$ 
\item[(ii)] $\forall j: \, X_j$ denotes the $j$´th column of $X.$ 
\item[(iii)] The column norm of $X$ is denoted $\|X\|_c$ and it is defined as $\|X\|_c :=  \underset{j}{\max}\|X_j\|_2.$ 
\item[(iv)] For a vector $\Lambda = ( \l_1, \dots ,\l_l)$ in $\bc^l,$ we let $\Lambda_{-}$ and $\Lambda_|$ denote respectively the $1\times l $ row matrix and the $l \times 1$ column matrix with entries from $\Lambda.$
\item[(v)] For any vector $\Lambda$ in $\bc^l$ the expression $\Delta(\Lambda)$ denotes the diagonal $l \times l$ matrix with entries from $\Lambda.$ 
\item[(vi)] The symbol $\bt$ denotes the torus, or unit circle,  in the complex plane.
\item[(vi)] The symbol $\ca_n$ denotes the $n-$dimensional abelian C*-algebra $C(\{1, \dots, n\}, \bc).$ 
\item[(vii)]
The complex $m \times n$ matrices are equipped with the  inner product given by $\langle X, Y \rangle := \mathrm{Tr}_n(Y^*X) .$ 
\end{itemize}
\end{definition}

We can now formulate our first factorization result.

 \begin{theorem} \label{ColBdFac} 
Let $ n$ be a  natural number,    then there exists a minimal positive real $c_n \leq \sqrt{k_G^\bc}$ such that for any natural number $m$ and  any $X$ in $M_{(m,n)}(\bc)$ there exist  a natural number $l \leq n^2 ,$ a matrix $C$ in $M_{(m,l)}(\bc) $ with $\|C\|_\infty \leq  c_n\|X\|_c$  a unit vector $\Lambda = ( \l_1, \dots , \l_l)$ in $\bc^l$ with positive entries and a matrix $Z \in M_{(l,m)}(\bt)$ such that $ X = C \Delta(\Lambda)Z.$ 

The supremum of the set $\{c_n: n \in \bn\} $ equals $\sqrt{k_G^\bc}.$  
\end{theorem}

\begin{proof}
Suppose that $\|X\|_c = 1,$ then by Proposition \ref{geo} there exists an operator $R$ in $\car_n$ such that $X^*X \leq k^\bc_G R.$ Since any element in $\car_n$ has the unit as its diagonal, the affine dimension of the convex set $\car_n$ is at most $2*(1/2)n(n-1)\leq n^2 -1,$ hence by Carathéodory's Theorem, $ R$ may be expressed as a convex combination of length $l$ with $l \leq n^2$ of the form \begin{equation} \label{ConvComb} 
R_{(i,j)} = \sum_{k=1}^l \l_k^2 \bar u^k_i u^k_j \text{ with }  \l_k > 0, \, \, \sum_{k=1}^l \l_k^2 =1,\, \,  u^k_j \in \bt 
\end{equation}
We can then define a matrix $Z$ in $M_{(l,n)}(\bt) $ by $Z_{(k,j)} = u^k_j,$  and it is easy to see that $R = (Z^*\Delta(\lambda)^*) (\Delta(\Lambda)Z)$ so the inequality $ X^*X \leq k^\bc_G R ,$ implies the existence of a matrix $C$ in $M_{(m,l)}(\bc)$ with $\|C\|_\infty \leq \sqrt{k_G^\bc}$ such that $X = C\Delta(\Lambda)Z,$ and the first part of the theorem has been proved.  

 Now let $h:= \sup_n\{c_n\},$ then $h \leq \sqrt{k_G^\bc}, $ and in order to prove the opposite inequality we return to the norm $\|.\|_T$ which we mentioned in the introduction right after equation (\ref{lG}). 
 
 We will now - based on the previous parts of this proof - show that for any $X$ in $M_{(m,n)}(\bc)$ we have $\|X\|_{\wedge(2, \infty)}  \leq h\|X\|_T,$ which will imply that $h \geq \sqrt{k_G^\bc}, $ because $\sqrt{k_G^\bc}$ is the minimal constant usable in equation (\ref{gT}). Then let $X$ be in $M_{(m,n)}(\bc)$ with $\|X\|_T = 1$ of rank $r,$   then  by item (iv) in Theorem 2.8 of \cite{C3}, or item (iv) in Theorem 3.2 of \cite{C4},  there exist a non negative unit vector $\Gamma$ in $\bc^m, $ a matrix $L$ in in $M_{(r,m)}(\bc)$ with $\|L\|_c = 1 $ and a matrix $R$ in $M_{(r,n)}(\bc)$  with $\|R\|_c = 1, $ such that $X= \Delta(\Gamma)L^*R. $ 
 We will apply Theorem   \ref{ColBdFac} to $R$ and use that $h \geq c_n,$  so $R$ has a decomposition $R = C \Delta(\Lambda) Z$ for a positive unit vector $\Lambda $ in $\bc^l,$ some $Z$ in $M_{(l,n)}(\bt)$ and a $C$ in $M_{(r,l)}(\bc)$  with $\|C\|_\infty \leq h.$  All the columns in $L$ have  length at most 1, so the 2-norm of $\Delta(\Gamma)L^* $ is at most $1,$ and for the matrix $M:= \Delta(\Gamma)L^*C$ we have $\|M\|_2 \leq h.$ Hence the sum of the squares of the 
 norms of the columns in 
$M$ satisfy $\sum_j\|M_j\|_2^2 \leq h^2. $ Let $N := \Delta(\Lambda)Z,$ then the sum of the squares of the  $\infty-$norms of the rows in $N$ satisfy $\sum_j \|_jN\|_\infty^2 = \sum_j \l_j^2 =  1.$ Then  if we consider $X=MN = \sum_j (M_j)_| (_jN)_-$ as an element in $ (\bc^m, \|.\|_2) \otimes ( \bc^n, \|.\|_\infty)$ equipped with the projective norm $\|.\|_{\wedge(  
2, \infty )},$ then we get by Cauchy-Schwarz' inequality \begin{equation} \notag
\|X\|_{\wedge(2, \infty)} \leq \sum_j \|M_j\|_2\|_jN\|_\infty  = \sum_j \|M_j\|_2\l_j \leq \|M\|_2 \sqrt{ \sum_j \l_j^2} \leq h.\end{equation} 
Since $X$ was chosen arbitrarily, except for the demand that $\|X\|_T = 1,$ it follows from (\ref{gT}) that $\sqrt{k_G^\bc} \leq h$ and the theorem follows.  \end{proof}

The statement, that $\sqrt{k_G^\bc} $ is the best possible constant working for any natural number $n,$ actually shows that the content of the Theorem \ref{ColBdFac} is an equivalent formulation of Grothendieck's little inequality. Further the proof offers a proof of the following corollary, which is a bit  stronger  than the statement $\|X\|_{\wedge(2,\infty)} \leq \sqrt{k_G^\bc}\|X\|_T.$ 

\begin{corollary}    
Let $X$ be a  complex $m \times n$ matrix, then there exist a  natural number $l \leq n^2$ an $l \times n$ matrix $Z$ with entries in $\bt$ and a complex $m \times l$ matrix $L$ such that $X=LZ $ and $\sum_t \|L_t\|_2 \leq \sqrt{k_G^\bc} \|X\|_T.$ 
\end{corollary} 
In our opinion, the potential for applications of  Theorem \ref{ColBdFac} lies in the fact, that the number  $n $ of columns only is used in the inequality $l \leq n^2,$ so if we do not care about such an estimate we may add as many extra columns as we want, and we can still get a matrix $C$ which works for all of the new columns, too. We have tried to look into this aspect, but so far in vain, except for Theorem \ref{ManyCol}, which you will find right after the following theorem. That next theorem is  a slightly weaker version of Theorem \ref{ColBdFac}, which has no upper bound on $l$ and also depends on an arbitrary positive real $\eps.$ On the other hand it has the advantage that the factor $\Delta(\Lambda)$ is a multiple of the identity.  This may be nice to have, we hope, in a possible  application of the theorem.  The following result also has some relations to the aspects of Grothendieck's inequalities, which deal with embeddings of classical  Banach spaces into each other \cite{Pi3}. We are not very familiar with these results, so we have not been able to deal seriously with this aspect.

\begin{theorem} \label{ColBd2} 
Let $X$ be in $M_{(m,n)}(\bc)$ with $\|X\|_c \leq 1$ and $0 < \eps < 1,$ then there exist a natural number $\hat l, $ a matrix $\hat C$ in $M_{(m,\hat l)}(\bc) $ with $\|\hat C\|_\infty < \sqrt{k^\bc_G}(1 + \eps)$ and a matrix $ \hat Z $ in $M_{(\hat l,n)}(\bt)$ such that $X = (\hat l)^{-1/2}\hat C \hat Z.$ 
\end{theorem}

\begin{proof}
Based on  Theorem \ref{ColBdFac} we find a natural number $l \leq n^2$ and  a factorization of $X$ as $X = C\Delta(\Lambda)Z,$ as described in the theorem.  
Let us fix a natural number $q $ which satisfies \begin{equation} \label{q}  q > 2 l \underset{1 \leq j \leq l}{\max}\{ \l_j^{-2}\} \eps^{-1}. \end{equation} Then determine natural numbers $p_j$ and non negative reals $\mu_j$ such that  \begin{equation} \label{lj} \forall j \in \{1, \dots,l\} : \quad \l_j^2 \, = \, \frac{p_j}{q} - \mu_j \text{ and } 0 \leq \mu_j < \frac{1}{q}.
\end{equation} The two previous equations imply that in most of the cases, $p_j$ will be  a large natural number,
\begin{equation} \label{pj} p_j = q\l_j^2 +q \mu_j > 2l\eps^{-1} > 2l.
\end{equation}  We can define a non negative real $\d_j$ by \begin{equation} \label{dj} \d_j := \sqrt{p_j} - \sqrt{p_j - q\mu_j},
\end{equation} 
For each index $j$ we have  $$p_j(p_j^{-1/2}q^{-1/2}) - \d_jq^{-1/2} = \sqrt{p_j/q} - \d_jq^{-1/2} = \sqrt{p_j/q -\mu_j } = \l_j,$$ so we  may express the $ 1 \times 1 $ complex matrix $(\l_j)$ as the matrix product of a $ 1 \times (p_j+1)$ complex row matrix and a  $(p_j+1) \times 1$ complex column matrix as shown below.

\begin{equation} \label{explj} 
(\l_j) = \begin{pmatrix}
p_j^{-1/2}  & \dots & p_j^{-1/2}& - \d_j
\end{pmatrix} \begin{pmatrix}
q^{-1/2} \\ \dots \\ q^{-1/2}\\ q^{-1/2}
\end{pmatrix}.
\end{equation}

Based on equation \ref{explj} we can now for each $j$  look at  the $j$'th column $(C\Delta(\Lambda))_j= \l_jC_j $ of $C\Delta(\Lambda)$  and replace this - single -  $m-$column by an $m \times (p_j+1)$ complex matrix $\tilde C_j $ such that the first $p_j$  columns in $\tilde C_j$ equals $p_j^{-1/2}C_j$ and the last one equals $-\d_jC_j.$

\begin{equation}
\tilde C_j = \begin{pmatrix} p_j^{-1/2} C_j & \dots & p_j^{-1/2} C_j& -\d_jC_j
\end{pmatrix}.
\end{equation} 
It follows from equation (\ref{explj}) that \begin{equation}
\l_j C_j = \tilde C_j \begin{pmatrix}
q^{-1/2}\\ \dots\\ q^{-1/2}\\q^{-1/2}
\end{pmatrix}.
\end{equation}

We can now start to construct the natural number $\hat l$ and the  matrices $\hat C $ and $\hat Z. $  After that we will prove that they have the properties promised. We define \begin{equation} \label{hatl} 
\hat l : =\sum_{j =1}^l (p_j +1) = l + \sum_{j = 1}^l p_j.
\end{equation}
The matrix $\hat C$ can not be constructed in one step, so we will first define an $m \times \hat l$ complex matrix $\tilde C$ which is obtained from $C\Delta(\Lambda)$ by replacing each column $ \l_jC_j$ with the $m \times (p_j +1)$ scalar matrix $\tilde C_j,$ or it may described as the concatenation of all the matrices $\tilde C_j,$  as shown in the next equation
\begin{equation} \label{tildeC}
\tilde C := \begin{pmatrix}
\tilde C_1| & \dots & |\tilde C_l
\end{pmatrix}.
\end{equation}  The matrix $\hat Z$ is obtained from $Z$ by replacing each row $_jZ$ by $p_j +1 $ copies of the very same row, so $\hat Z $ is a $\hat l \times n $ matrix with entries in $\bt,$ and it looks like 
\begin{equation}
\hat Z = \begin{pmatrix}
_1Z \\
\dots\\ _1Z\\ _2Z \\ \dots \\ _lZ \\ \dots \\ _lZ
\end{pmatrix}
\end{equation}
It now follows from the constructions above that 
\begin{equation} \label{tildehat}
X = q^{-1/2} \tilde C \hat Z,
\end{equation} and we must define \begin{equation} \label{hatC} 
\hat C := \hat l^{1/2}q^{-1/2} \tilde C \end{equation} in order to complete the proof, we then have to  show that the matrix $\hat C$ satisfies $\|\hat C\|_\infty \leq \sqrt{k^\bc_G}(1 + \eps). $ 

We start by estimating the norm of $\tilde C \tilde C^*,$   and to that end we let $(C_j)^*$ denote row matrix which is adjoint of the column matrix $C_j.$  Then \begin{align} \label{TildeCTildeC*} 
\tilde C \tilde C^*& = \sum_{j=1 }^l \big( p_j (1/p_j) + \d_j^2\big) C_j(C_j)^* \leq (1 + \underset{j}{\max}\{\d_j^2\} )\sum_{i=1}^l C_i(C_i)^* \\ \notag &= (1 + \underset{j}{\max}\{\d_j^2\} )CC^* 
\end{align} 
To estimate $\d_j^2$ we return to the equations (\ref{dj})  and  (\ref{lj}), from where we get  we get for each index $j$ 
\begin{align} \label{dj<e2} 
\d_j^2 &\leq (\sqrt{p_j} - \sqrt{p_j -1})^2 = (\frac{1}{\sqrt{p_j} + \sqrt{p_j -1}})^2 \\ \notag & < \frac{1}{p_j}  \leq  \frac{1}{2l} \eps \leq \frac{1}{2} \eps.
\end{align} For any real $x \geq -1 $ we have $\sqrt{1+x} \leq  1 + (1/2)x,$ then by (\ref{TildeCTildeC*}) we obtain \begin{equation} \label{NormTildeC} 
\|\tilde C\|_\infty \leq \sqrt{k^\bc_G}\sqrt{1 + (1/2)\eps} \leq  \sqrt{k^\bc_G}(1 + (1/4)\eps). 
\end{equation}

To estimate $\sqrt{\hat l} / \sqrt{q}$ we return to the definition of $\hat l$ from (\ref{hatl}) and the equation (\ref{lj}). Then we recall that $\sum_j \l_j^2 = 1, $ so we get \begin{align} \notag  q = q \sum_j \l_j^2 &= \sum_j(p_j - q\mu_j ) = (\hat l -l)  - \sum_j q\mu_j,\\ \notag
\text{ and by (\ref{lj}), } 0 & \leq \sum_j q\mu_j \leq  l \\ \notag \text{ so } \hat l  - 2l & \leq q \leq \hat l - l
\\ \text{ and }        q + l & \leq  \hat l \leq q + 2l . \label{l/q} \end{align}
 Then we get from (\ref{l/q}) and  (\ref{q}) that 
 \begin{align} \notag \sqrt{\hat l/ q} & \leq (1+ 2l /q)^{(1/2)} \leq 1 + l/q \\ \notag & \leq 1 + \frac{l}{2l \max\{\l_j^{-2}\} \eps^{-1}} \leq 1+ \frac{1}{2} \eps, 
\end{align}
and then  for $\hat C := \sqrt{\hat l/q} \, \tilde C$ we get, since $0 < \eps < 1$   $$\|\hat C\|_\infty \leq \sqrt{k^\bc_G}(1+ (1/4) \eps) (1+ (1/2)\eps) < \sqrt{k^\bc_G}(1+  \eps),$$ and the theorem follows.   
\end{proof}

Right before this theorem we mentioned the possibility of extending an $m \times n$ matrix by an $m \times p$ matrix $Y$ via concatenation, $( X|Y) ,$ and then use the factorization theorem on the extended matrix. The next theorem studies the effect of this program in the case, where the finite number  of columns in $Y$ forms an $\eps-$dense subset of  the unit sphere in  $\bc^m.$

\begin{theorem} \label{ManyCol} 
Let $X$ be a complex $m \times n$ matrix with $\|X\|_c \leq 1, $ and $0 < \eps <1,$ then there exist a natural number $l,$ a  matrix $C$ in $M_{(m,l)}(\bc)$ with $$(1- \eps )^2 I_m \leq CC^* \leq k_G^\bc(1 + \eps)^2I_m$$
 a matrix $Z^X$ in $M_{(l, n)}(\bt)$ such that $X = l^{-(1/2)} CZ^X$ and a matrix $Z^{I_m} $ in $M_{(l,m)} (\bt)$ such that $ I_m = l^{-(1/2)}CZ^{I_m}.$ 
\end{theorem}

\begin{proof}
Let $\cas := \{\s_1, \dots, \s_p\}$ denote a set of unit vectors in $\bc^m$ which is $\eps-$dense in the unit sphere of $\bc^m.$ Then we let $S$ denote the  complex  the $m \times p $ matrix  which is defined by the equations $S_j := (\s_j)_|.$ We can then define a complex $m \times(n +m +p)$ matrix by concatenating the 3 matrices $X, I_m$ and $S$ as $Y := (\, X \,| \,I_m \,|
 \, S\, ), $ and in this way we have obtained a matrix which satisfies $\|Y\|_c \leq 1.$  When we apply Theorem \ref{ColBd2}
  to $Y,$ and we get a natural number $l,$ a complex $m \times l$  matrix $C$ with $\|C\|_\infty \leq \sqrt{k_G^\bc}(1 + \eps)$ 
 and  matrices $Z^X $ in $M_{(l,n)}(\bt),$ $Z^{I_m}  $ in $M_{(l,m)}(\bt),$ and  $Z^S $ in $M_{(l,p)}(\bt),$ such that \begin{equation} \label{C} (i):\, X= l^{-(1/2)}CZ^X,\,\,  (ii):\, I_m = l^{-(1/2)}CZ^{I_m},\,\,  (iii):\, S = l^{-(1/2)}CZ^S.\end{equation}
 Now let $VA  = C^*$ be the polar decomposition of $C^*,$ then we see from equation  (\ref{C}) (ii), that  the range of $C$ is all of $\bc^m$, so the positive matrix $A$ in $M_m(\bc)$ is invertible and we can define a matrix $D$ in $M_{(l, m)}(\bc)$ by $D:=VA^{-1}. $ 
 For each $\s_j$ we get by (\ref{C}) (iii) that $$(D\s_j)_| = DS_j = l^{-(1/2)}DCZ^S_j = VV^* l^{-(1/2)}Z^S_j,$$ and then $\|D\s_j\|_2 \leq 1,$ since $\|l^{-(1/2)}Z^S_j\|_2 =1.$
 
 The vectors $\s_j$ are $\eps-$ dense in the unit sphere  of $\bc^m,$ so by elementary techniques,  we get that  any unit vector $ \s$ in $\bc^m$  may be written as a norm convergent sum $$\s = \s_{j_0} + \sum_{q = 1 }^\infty t_q\s_{j_q }, \text{ with } 0 \leq t_q \leq \eps^q,\, \s_q \in \cas$$  and then $\|D\s\|_2 \leq 1/(1-\eps).$ This inequality and the inequality $\|C\|_\infty \leq \sqrt{k_G^\bc}(1+\eps)$ then gives \begin{equation} \label{CC*}
 (1-\eps)^2I_m \leq CC^* \leq k_G^\bc(1+\eps)^2I_m,
 \end{equation} and the theorem follows.
 
 \end{proof}

\section{Proofs  of Grothendieck's inequalities based on  Theorem \ref{ColBdFac}}
 
In Theorem 3.2 of \cite{C4} and later on in equation (4.1) of \cite{C1} we remarked that
\begin{align} &\text{The constant } k_G^\bc = \frac{\pi}{4} \text{ is the minimal positive real such that }\\ \notag &\forall n \in \bn\, \forall P \in M_n(\bc ) \text{ with } P\geq 0:  \quad \|P\|_{cbB} \leq k_G^\bc\|P\|_B. \end{align}  
This inequality is an equivalent formulation of Grothendieck's
{ \em little inequality} based on the fact that the operator norm satisfies $\|X^*X\|_\infty = \|X\|_\infty^2.$   When examining this observation  we showed in \cite{C1}, that this implies a general validity of the big inequality, with a constant dominated by $k_G^\bc/(2- k_g^\bc).$  Here we show that both of these results also may be deduced from the content of Theorem \ref{ColBdFac}.  
 
 \begin{proposition}
Let $X$ in $M_n(\bc)$ be a positive matrix, \newline  then $\|X\|_{cbB} \leq k_G^\bc \|X\|_B.$ \newline
Let $X$ be in $M_{(m,n)}(\bc),$ then $\|X\|_{cbB} \leq  
 k_G^\bc/(2- k_G^\bc)\|X\|_B.$
 \end{proposition} 
 
 \begin{proof}
 
 Given a positive $X$ in $M_n(\bc),$ we may and will assume that $\|X\|_{cbB} =1.$  
 By equation (3.5) of \cite{C3} we know that the equation $\|X\|_{cbB} =1,$  implies that there exists a positive $n \times n$ matrix $Y$ of Schur multiplier norm 1 such that Tr$_n(YX) =1. $   Since $\|Y\|_S = 1$ and $Y \geq 0$ we know from Schur's pioneering work \cite{Sc}, that $0 \leq \mathrm{diag}(Y) \leq I_n. $ 
 The square root $Y^{(1/2)} $ is then column bounded by 1, and by Theorem \ref{ColBdFac}, it  has a factorization $$ Y^{(1/2)} = C \Delta(\Lambda)Z$$ with  the properties that $\l_i> 0 ,$ $\sum\l_i^2 = 1,$  $\|C\|_{\infty}^2 \leq k_G^\bc$  and the entries of $Z$ are in $\bt.$   Based on this we can estimate as follows  
 \begin{align*}
 \|X\|_{cbB} =1 &= \mathrm{Tr}_n( XY) = \mathrm{Tr}_n(XZ^*\Delta(\Lambda)C^*C \Delta(\Lambda) Z) \\
 &= \mathrm{Tr}_l(\Delta(\Lambda) ZXZ^*\Delta(\Lambda)C^*C ) \\
 & \leq k_G^\bc \mathrm{Tr}_l(\Delta(\Lambda)Z X Z^*\Delta(\Lambda))\\&= k_G^\bc \sum_{i=1}^l \lambda_i^2B_X(_iZ, \overline{(_iZ)}) \\ &\leq k_G^\bc \|X\|_B,
 \end{align*}
and the first statement is proved. It is worth to remark that  for a general complex $m \times n$ matrix $Y$ the matrix $P:= Y^*Y$ is a positive matrix, and the statement just proven then gives that $\|Y\|_{cbF} \leq \sqrt{k_G^\bc}\|Y\|_F.$  

Now let $X$ in $M_{(m,n)} (\bc) $ satisfy $\|X\|_{cbB} =1, $ then by Proposition 3.2 of \cite{C3} there exist a non negative unit vector $\eta$ in $\bc^m$ a matrix $B$ in $M_{(m,n)}(\bc)$ of operator norm 1 and a non negative unit vector $\xi $ in $\bc^n$ such that $X= \Delta(\eta)B \Delta(\xi).$  The duality described in equation (3.5) of \cite{C3} shows that there exists a  matrix $Y$ in $M_{(m,n)}(\bc)$ such that $\|Y\|_S = 1$ and Tr$_n(Y^*X) = 1.$ By Theorem 2.7 of \cite{C3} there exists a natural number $r,$ and matrices $L $ in $M_{(r,m)}(\bc) $ with $\|L\|_c =1 $ and $R$ in $M_{(r,n)}(\bc) $ with $\|R\|_c = 1$ such that $Y = L^*R.$
 We can then construct an $r  \times (m+n)$ matrix $M$ by concatenating $L$ and $R$ into $ M:= (L | R).$ 
 This matrix $M$ has  column norm equal to 1, and an application of Theorem \ref{ColBdFac} gives a factorization 
 $$ (L |R ) = M = C \Delta(\Lambda)(Z^L | Z^R), \,  Z^L \in M_{(l,m)}(\bt) \text{ and } Z^R \in M_{(l, n)}(\bt).$$ We can then define 2 matrices $S$ and $T$  with
  \begin{align*}S \in  M_{(r,n)}(\bc) &\text{ given by } S :=  L\Delta(\eta) B =  C\Delta(\Lambda)Z^L \Delta(\eta) B, \\
   T \in M_{(r,n)}(\bc) &\text{ given by }  T := R\Delta(\xi) =  C\Delta(\Lambda)Z^R \Delta(\xi) .
   \end{align*} 
 Since $\|L\|_c = 1,$ $\|\eta\|_2 = 1$ and $\|B\|_\infty = 1 $ we get from the first equality sign above that $\|S\|_2 \leq 1.$ 
 Similarly we get $\|T\|_2 \leq 1,$ and by computation we find $\langle S, T \rangle = \mathrm{Tr}_n(T^*S) = \mathrm{Tr}_n (Y^*X) =1. $ By Cauchy Schwarz' inequality we get that $\|S\|_2 = \|T\|_2 = 1$  and $S=T. $  
Let us define two matrices $\hat S $ and $\hat T$ in $M_{(l,n)}(\bc) $ by $$ \hat S := \Delta(\Lambda)Z^L \Delta(\eta) B \text { and } \hat T := \Delta(\Lambda) Z^R\Delta(\xi).$$  Finally let  $Q$ in   $M_l(\bc)$ denote  the matrix of the {\em orthogonal } support projection of $C.$  Since the entries of $Z^L$ and $Z^R$ are in $\bt$ and the vectors $\Lambda, \eta$ and $\xi$ are unit vectors we have $\|\hat S\|_2 \leq 1 $ and $\| \hat T\|_2 = 1.$  On the other hand the previous equation gives
\begin{align*} 
&Q \hat S = Q \hat T, \text{ and } \|Q\hat S \|_2 =  \|Q\hat T\|_2 \geq (k_G^\bc)^{-(1/2)} \text{ by } \|C\|_\infty \leq \sqrt{k_G^\bc}.\\ \notag 
& \text{ Since } 1 \geq \|\hat S\|_2^2 = \|Q \hat S\|_2^2 + \|(I_l - Q) \hat S\|_2^2 \text{ and similarly for } \hat T\\ 
&\|(I_l - Q) \hat S \|_2 \leq (1 - (k_G^\bc)^{-1})^{(1/2)}, \,  \|(I_l - Q) \hat T\|_2 \leq (1 - (k_G^\bc)^{-1})^{(1/2)}.
   \end{align*} 
These equations imply that 
\begin{equation} \label{hatS}
|\langle \hat S , \hat T\rangle | \geq (k_G^\bc)^{-1} - (1- (k_G^\bc)^{-1})  = 2 (k_G^\bc)^{-1} -1.
\end{equation}
 On the other hand \begin{align}  \label{BX} 
 |\langle \hat S, \hat T \rangle | &=| \mathrm{Tr}_n (\Delta(\xi) (Z^R)^* \Delta(\Lambda)^2 Z^L \Delta(\eta) B)| \\ \notag &= |\mathrm{Tr}_l ( \Delta(\Lambda) Z^L (\Delta(\eta)B\Delta(\xi)) (Z^R)^* \Delta(\Lambda))\\ \notag &= |\mathrm{Tr}_l ( \Delta(\Lambda) Z^L X (Z^R)^* \Delta(\Lambda)) |\\ \notag  &= | \sum_{i =1}^l \l_i^2 B_X (_iZ_L , \overline{_iZ_R})| \\ \notag  &\leq \|X\|_B. 
 \end{align}
 Hence by (\ref{hatS}) and (\ref{BX}) we get
 $\|X\|_B \geq 2(k_G^\bc)^{-1} - 1$ and then $K_G^\bc \leq (2(k_G^\bc)^{-1} - 1)^{-1} = k_G^\bc/(2 - k_G^\bc),$ and the proposition follows. 
 \end{proof}
 Numerically the proposition shows that $1.273 < 4/\pi =  k_G^\bc \leq  K_G^\bc < 1.752,$ and  we know from \cite{Ha3} and \cite{Pi3} Section 4 that $  k_G^\bc < 1.274  <  1.338 < K_G^\bc < 1.405.$

\end{document}